\documentclass[preprint,12pt]{elsarticle}

\usepackage{amsmath,amssymb,amsthm,mathtools}
\usepackage{graphicx}
\usepackage{tikz}
\usepackage{subcaption}
\usepackage{placeins}
\usepackage{enumitem,xcolor,microtype}
\usepackage[margin=2.5cm]{geometry}
\usepackage[colorlinks=true,linkcolor=blue!55!black,urlcolor=blue!55!black]{hyperref}
\usepackage{xurl}

\newtheorem{theorem}{Theorem}[section]
\newtheorem{lemma}[theorem]{Lemma}
\newtheorem{proposition}[theorem]{Proposition}
\newtheorem{corollary}[theorem]{Corollary}
\theoremstyle{definition}
\newtheorem{definition}[theorem]{Definition}
\newtheorem{remark}[theorem]{Remark}
\newtheorem{example}[theorem]{Example}
\numberwithin{equation}{section}

\newcommand{\sig}{\sigma}
\newcommand{\ph}{\varphi}
\newcommand{\bbone}{\mathbf1}
\newcommand{\scoronal}{\eta}
\newcommand{\PSD}{\succeq0}
\newcommand{\Sp}{\operatorname{Sp}}
\newcommand{\Spstar}{\operatorname{Sp}^{*}}
\newcommand{\join}{\vee}
\newcommand{\udot}{\mathbin{\dot\cup}}

\begin{document}

\begin{frontmatter}

\title{A complete characterization of graphs with second largest eigenvalue
at most \texorpdfstring{$(\sqrt5-1)/2$}{(sqrt(5)-1)/2}}

\author{Chaohui Chen\texorpdfstring{\textsuperscript{a}}{a}}
\ead{cch956833210@163.com}
\author{Hangxi Cha\texorpdfstring{\textsuperscript{a}}{a}}
\ead{2433949@tongji.edu.cn}
\author{Haiying Shan\texorpdfstring{\textsuperscript{a}}{a}}
\ead{shan_haiying@tongji.edu.cn}
\author{Muhuo Liu\texorpdfstring{\textsuperscript{b}}{b}}
\ead{liumuhuo@163.com}
\affiliation{organization={
School of Mathematical Sciences, Key Laboratory of Intelligent Computing and Applications\\
(Ministry of Education), Tongji University},
            city={Shanghai},
            postcode={200092},
            country={P.R. China}}
\affiliation{organization={Department of Mathematics, South China Agricultural University},
            city={Guangzhou},
            country={P.R. China}}
\nonumnote{This work is partially supported by the National Natural Science Foundation of China.}

\begin{abstract}
In 1993, Cao and Hong posed the problem of
characterizing graphs whose second largest adjacency eigenvalue is
greater than $1/3$ and less than $(\sqrt5-1)/2$.  Subsequent work extended
the question to the non-strict bound.  We give a complete structural characterization
of connected graphs $G$ with $|V(G)|\ge2$ satisfying
$\lambda_2(G)\le(\sqrt5-1)/2$.
Their complements are described by disjoint unions of graphs from
an explicit family with least adjacency eigenvalue at least
$-(1+\sqrt5)/2$, possibly together with a cone or a thin or thick
spider constructed from disjoint unions of members of the same
family, subject to explicit numerical inequalities.
This provides a positive characterization of the class, going
beyond the classical result of Cvetkovi\'c and Simi\'c that the bound
$\lambda_2\le(\sqrt5-1)/2$ admits only finitely many minimal
forbidden induced subgraphs.
\end{abstract}


\end{frontmatter}

\section{Introduction}\label{sec:intro}
Let $G$ be a finite simple graph of order $n$, with adjacency eigenvalues
$\lambda_1(G)\ge\lambda_2(G)\ge\cdots\ge\lambda_n(G)$.
Determining graphs with a prescribed upper bound on $\lambda_2(G)$
is a long-standing problem in spectral graph theory.
The case $\lambda_2(G)\le1$ was studied by Cvetkovi\'c
\cite{Cvetkovic1982}, and several restricted classes have since been
characterized; see, for example,
\cite{LiuChenStanic2021,LiuChenStanic2024,ChenShanLiuStanic2026,Cheng2026}.
For general background, see
\cite{CvetkovicRowlinsonSimic2010,Stanic2015}.

Let $\sig=(\sqrt5-1)/2$ denote the reciprocal of the golden ratio,
for which $\lambda_2(P_4)=\sig$.  Cao and Hong
\cite{CaoHong1993} posed the problem of characterizing graphs with
$1/3<\lambda_2(G)<\sig$.  Cvetkovi\'c and Simi\'c subsequently
investigated graphs satisfying $\lambda_2(G)\le\sig$ by means of
hereditary properties and forbidden induced subgraphs
\cite{CvetkovicSimic1995}.  Simi\'c established finiteness of the
forbidden induced-subgraph family for the strict bound
\cite{Simic1995}; further results on forbidden induced subgraphs appeared in
\cite{CvetkovicSimic2000}.  The case $\lambda_2(G)<1/2$ was
characterized by Wu, Qian and Peng \cite{WuQianPeng2023}.
More recently, Liu, Chen, Stani\'c and Shan
\cite{LiuChenStanicShan2025} characterized graphs with
$\lambda_2(G)\le\sig$ when the clique number is at least $54$, or at
least $8$ in the presence of a pendant vertex.

We remove the clique-number restriction by giving a complete structural
characterization of connected graphs $G$ of order at least $2$ with
$\lambda_2(G)\le\sig$.  Our characterization describes their complements
in terms of disjoint unions, cones, and thin and thick spiders constructed
from an explicit family of graphs, subject to explicit numerical
inequalities.

We work with the complement $H=\overline G$ and put
$\ph=1+\sig=1/\sig$.  We first determine the family $\mathcal H_+$
of connected graphs with least adjacency eigenvalue at least $-\ph$.
This classifies the connected components of $H$ whenever
$\lambda_{\min}(H)\ge-\ph$.  For the remaining case, we identify the
five minimal forbidden induced subgraphs for this least-eigenvalue
bound whose complements satisfy $\lambda_2\le\sig$.
Neighbourhood analysis around these graphs yields the cone and spider
decompositions needed for the characterization.

Finally, a rank-one perturbation argument, together with shifted-coronal
formulas for these constructions, gives necessary and sufficient
numerical conditions for the bound $\lambda_2(G)\le\sig$.
These conditions can be checked from the component parameters,
without computing the adjacency eigenvalues of $G$.

Section~\ref{sec:matrix} gives the spectral preliminaries, introduces
the graph families and their shifted coronals, and states the main theorem.
Section~\ref{sec:positive-seeds} characterizes the connected graphs with
least eigenvalue at least $-\ph$ and identifies the five minimal forbidden
induced subgraphs whose complements satisfy $\lambda_2\le\sig$.
Sections~\ref{sec:a1} and~\ref{sec:clawfree} establish the cone and spider
decompositions according as the complement contains a claw or is claw-free.
Section~\ref{sec:global} derives the inertia and coronal formulas and
completes the proof of the main theorem.

\section{Preliminaries and the main theorem}\label{sec:matrix}
All graphs are finite and simple.  For a graph $G$, write $V(G)$,
$E(G)$, and $A(G)$ for its vertex set, edge set, and adjacency matrix,
and $|V(G)|$ for its order.  For $U\subseteq V(G)$, let $G[U]$ be
its induced subgraph, and put $G-v=G[V(G)\setminus\{v\}]$.
For disjoint $X,Y\subseteq V(G)$, $G[X,Y]$ denotes the bipartite graph
with parts $X,Y$ consisting of the edges of $G$ between them.
Write $N_G(v)=\{u\in V(G):uv\in E(G)\}$ and
$N_G[v]=N_G(v)\cup\{v\}$ for the open and closed neighbourhoods,
and put $N_G(U)=\bigcup_{v\in U}N_G(v)$.

Unless otherwise stated, $G$ is a connected graph of order at least $2$,
$H=\overline G$, and edges and neighbourhoods refer to $H$.
A \emph{non-neighbour} of $v$ is a vertex $u\ne v$ with $uv\notin E(H)$.

We write $K_n$, $C_n$, and $P_n$ for the complete graph, cycle, and path
on $n$ vertices.  The graph $K_{1,3}$ is the \emph{claw}, and a graph
is \emph{claw-free} if it contains no induced claw.
A \emph{leaf}, or \emph{pendant vertex}, has degree one, and a vertex is \emph{universal} if it is
adjacent to every other vertex of the graph.
For graphs on disjoint vertex sets, $H_1\udot H_2$ and $H_1\vee H_2$
denote their disjoint union and join, respectively; the join is obtained
by adding all edges between the two vertex sets.  We write $kG$ for the disjoint
union of $k$ copies of $G$, and $\bigsqcup$ for indexed disjoint unions.
Empty disjoint unions and empty sums are allowed, with an empty sum
equal to zero.  The symbol $\varnothing$ denotes either the empty set
or the empty graph, as appropriate.
Disjoint vertex sets are \emph{complete} (respectively,
\emph{anticomplete}) to each other if all (respectively, no) edges
between them are present.

For nonnegative integers $a,b$, let $J_{a,b}$ be the $a\times b$
all-ones matrix, with the empty-block convention when $ab=0$,
and put $J_a=J_{a,a}$.  Write $I_a$ for the identity matrix of order
$a$ and $\bbone_a$ for the all-ones vector of length $a$.
Subscripts are omitted when dimensions are clear;
$J=\bbone\bbone^{\mathsf T}$.
For a real symmetric matrix $M$, we write $M\succeq0$ if $M$ is
positive semidefinite, and let $M^\dagger$ denote its Moore--Penrose
inverse~\cite{BenIsraelGreville2003}.  If
$M=U\operatorname{diag}(\mu_1,\ldots,\mu_n)U^{\mathsf T}$ is an orthogonal
diagonalization, then
\[
M^\dagger=U\operatorname{diag}(\mu_1^\dagger,\ldots,\mu_n^\dagger)U^{\mathsf T},
\qquad
\mu_i^\dagger=
\begin{cases}
1/\mu_i,&\mu_i\ne0,\\
0,&\mu_i=0.
\end{cases}
\]
The space $\operatorname{range}(M)$ is the column space of $M$.
The \emph{inertia} of $M$ is
$\operatorname{In}(M)=(n_+(M),n_0(M),n_-(M))$, the numbers of positive,
zero and negative eigenvalues, counted with multiplicity.

For a nonempty graph $G$, let $\lambda_{\min}(G)$ denote its least
adjacency eigenvalue.  We call $A(G)+\ph I$ the \emph{shifted adjacency
matrix} of $G$.  Thus $A(G)+\ph I\PSD$ if and only if
$\lambda_{\min}(G)\ge-\ph$.

\begin{lemma}[\cite{Haemers1995}]\label{lem:induced-spectral-bounds}
Let $F$ be a graph and $F'$ a nonempty induced subgraph of $F$.  Then
$\lambda_{\min}(F')\ge\lambda_{\min}(F)$.  If $|V(F')|\ge2$, then also
$\lambda_2(F')\le\lambda_2(F)$.  Consequently, both bounds
$\lambda_{\min}\ge-\ph$ and $\lambda_2\le\sig$ are inherited by induced
subgraphs for which the relevant eigenvalue is defined.
\end{lemma}

\begin{lemma}[\cite{HornJohnson2013}]\label{lem:weyl}
For real symmetric matrices $A,B$ of order at least $2$, with
eigenvalues ordered nonincreasingly,
\[
\lambda_2(A+B)\le\lambda_2(A)+\lambda_1(B).
\]
\end{lemma}

\begin{lemma}[\cite{CarlsonHaynsworthMarkham1974,Haynsworth1968}]\label{lem:schur-congruence}
Let $A,D$ be real symmetric matrices and let $B$ satisfy
$\operatorname{range}(B)\subseteq\operatorname{range}(A)$.
Then
\[
\operatorname{In}\!\begin{pmatrix}A&B\\B^{\mathsf T}&D\end{pmatrix}
=
\operatorname{In}\!\begin{pmatrix}
A&0\\0&D-B^{\mathsf T}A^\dagger B
\end{pmatrix}.
\]
Here $D-B^{\mathsf T}A^\dagger B$ is the \emph{generalized Schur complement}
with respect to $A$; when $A$ is invertible, this is the usual Schur
complement $D-B^{\mathsf T}A^{-1}B$.
In particular,
\[
n_-\!\begin{pmatrix}A&B\\B^{\mathsf T}&D\end{pmatrix}
=n_-(A)+n_-(D-B^{\mathsf T}A^\dagger B).
\]
\end{lemma}

\begin{figure}[!htbp]
\centering
\resizebox{0.52\textwidth}{!}{%
\begin{tikzpicture}[
  line width=0.8pt,
  v/.style={circle,draw=black,fill=black,inner sep=0pt,minimum size=2mm},
  lab/.style={font=\small}
]
% ---------- Sp(P_3,3) : thin spider ----------
\begin{scope}[shift={(0,0)}]
% Vertex positions for the thin spider
% T (left arc)
\node[v] (t1) at (-0.7,1.212) {};
\node[v] (t2) at (-1.4,0) {};
\node[v] (t3) at (-0.7,-1.212) {};
% L = P_3 (vertical column on the right)
\node[v] (l1) at (1.4,1.212) {};
\node[v] (l2) at (1.4,0) {};
\node[v] (l3) at (1.4,-1.212) {};
% Q: degree-one vertices (vertical column on the left)
\node[v] (q1) at (-2.3,1.212) {};
\node[v] (q2) at (-2.3,0) {};
\node[v] (q3) at (-2.3,-1.212) {};
\draw (l1)--(l2)--(l3);
\draw (t1)--(t2)--(t3)--(t1);
% L complete to T
\foreach \l in {l1,l2,l3}{
  \foreach \t in {t1,t2,t3}{\draw (\l)--(\t);}
}
% Q-T perfect matching
\draw (q1)--(t1) (q2)--(t2) (q3)--(t3);
\end{scope}
% ---------- Sp*(P_3,3) : thick spider ----------
\begin{scope}[shift={(4.2,0)}]
% T (right arc)
\node[v] (t1b) at (0.7,1.212) {};
\node[v] (t2b) at (1.4,0) {};
\node[v] (t3b) at (0.7,-1.212) {};
% L = P_3 (vertical column on the left)
\node[v] (l1b) at (-1.4,1.212) {};
\node[v] (l2b) at (-1.4,0) {};
\node[v] (l3b) at (-1.4,-1.212) {};
% Q: degree-two vertices (vertical column on the right)
\node[v] (q1b) at (2.3,-1.212) {};
\node[v] (q2b) at (2.3,0) {};
\node[v] (q3b) at (2.3,1.212) {};
\draw (l1b)--(l2b)--(l3b);
\draw (t1b)--(t2b)--(t3b)--(t1b);
% L complete to T
\foreach \l in {l1b,l2b,l3b}{
  \foreach \t in {t1b,t2b,t3b}{\draw (\l)--(\t);}
}
% complement of the matching: q_i misses only its matched t_i
\draw (q1b)--(t2b) (q1b)--(t3b)
      (q2b)--(t1b) (q2b)--(t3b)
      (q3b)--(t1b) (q3b)--(t2b);
\end{scope}
\end{tikzpicture}}
\caption{The thin spider $\Sp(P_3,3)$ (left) and the thick spider
$\Spstar(P_3,3)$ (right).}
\label{fig:spiders}
\end{figure}
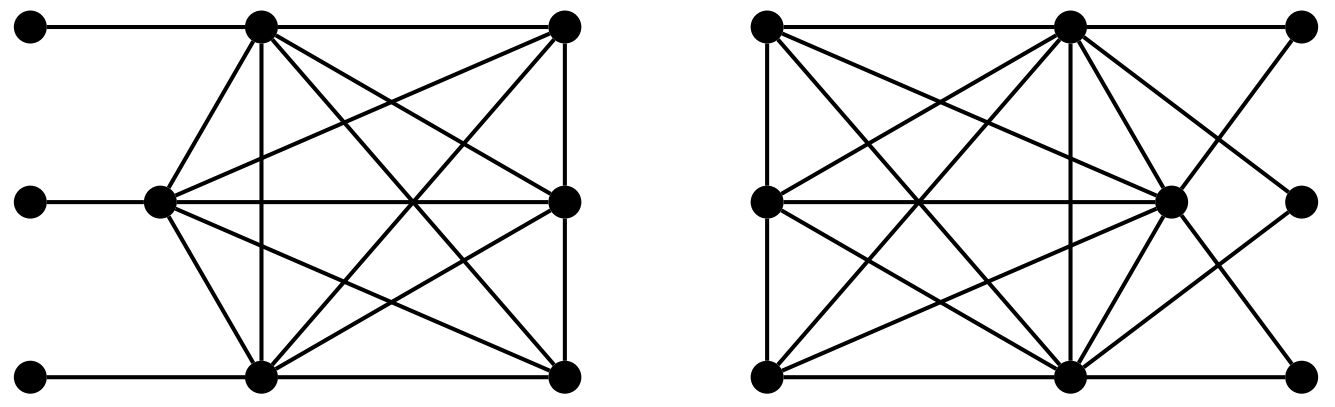

\begin{definition}\label{def:spider-families}
For a graph $L$, the graph $K_1\join L$ is the \emph{cone} over $L$,
with \emph{apex} its unique vertex outside $L$.
For an integer $t\ge1$, an \emph{$(L,t)$-spider} is a graph $G$ with a partition
$V(G)=V(L)\cup T\cup Q$ such that $G[V(L)]=L$, $|T|=|Q|=t$, $T$ is a
clique, $Q$ is independent, $V(L)$ is complete to $T$ and anticomplete
to $Q$, and $G[T,Q]\cong tK_2$ or $\overline G[T,Q]\cong tK_2$.
The two cases are called \emph{thin} and \emph{thick spiders}, denoted
by $\Sp(L,t)$ and $\Spstar(L,t)$, respectively; see
\cite{JamisonOlariu1992,BrandstadtDraganLeMosca2002} and
Figure~\ref{fig:spiders}.
\end{definition}

For $t=1$,
$\Spstar(L,1)\cong(K_1\join L)\udot K_1$ and
$\Sp(L,1)\cong K_1\join(L\udot K_1)$; for $t=2$,
$\Sp(L,2)\cong\Spstar(L,2)$.

Write $T=\{u_1,\ldots,u_t\}$ and $Q=\{v_1,\ldots,v_t\}$ so that
$u_iv_i\in E(\Sp(L,t))$ for each $i$.  Put $n=|V(L)|$.
With rows and columns indexed by $V(L),T,Q$ in this order,
\[
A(\Sp(L,t))=
\begin{pmatrix}
A(L)&J_{n,t}&0\\
J_{t,n}&J_t-I_t&I_t\\
0&I_t&0
\end{pmatrix}.
\]
For the thick spider, label the vertices so that $u_iv_i$ are the
missing edges between $T$ and $Q$.  Its adjacency matrix is obtained
by replacing the two off-diagonal blocks $I_t$ above by $J_t-I_t$.

\FloatBarrier

For integers $t\ge1$ and $0\le j\le t$, let $W_{t,j}$ be obtained from
$K_t$ by attaching $j$ leaves to distinct vertices.  In particular,
$W_{t,0}=K_t$.

For a graph $H$, the \emph{coronal} introduced by McLeman and McNicholas
\cite{McLemanMcNicholas2011} is
\[
\Gamma_H(\lambda)=\bbone^{\mathsf T}(\lambda I-A(H))^{-1}\bbone,
\]
whenever $\lambda I-A(H)$ is invertible.  We use the following
\emph{shifted coronal} at $-\ph$:
\[
\scoronal(H):=\bbone^{\mathsf T}(A(H)+\ph I)^\dagger\bbone.
\]
If $A(H)+\ph I$ is invertible, then $\scoronal(H)=-\Gamma_H(-\ph)$.
The Moore--Penrose inverse makes the definition valid also when
$A(H)+\ph I$ is singular.  If
$\bbone\in\operatorname{range}(A(H)+\ph I)$, then
$\scoronal(H)=\bbone^{\mathsf T}x$ for any solution of
$(A(H)+\ph I)x=\bbone$.  The value is independent of the chosen
solution, since $\bbone$ is orthogonal to the kernel of the symmetric
matrix $A(H)+\ph I$.

The shifted coronal is additive over disjoint unions: for any graphs
$L_1,\ldots,L_m$ with $m\ge0$,
\[
\scoronal\!\left(\bigsqcup_{j=1}^m L_j\right)
=\sum_{j=1}^m\scoronal(L_j),
\]
since $A(\bigsqcup_{j=1}^m L_j)+\ph I$ and its Moore--Penrose inverse
are block diagonal.  In particular, $\scoronal(\varnothing)=0$.

Put $\mathcal W=\{W_{t,j}:t\ge1,\ 0\le j\le t\}$ and define
\[
\mathcal H_+=\mathcal W\cup\{K_1\join P_3,K_1\join P_4,K_1\join C_5,C_5,
\overline{W_{3,3}},K_1\join(K_2\udot K_2),K_1\join(K_2\udot K_3)\}.
\]
\begin{lemma}\label{lem:positive-kernel-coronal}
For every $H\in\mathcal H_+$, we have $\lambda_{\min}(H)\ge-\ph$ and
$\bbone\in\operatorname{range}(A(H)+\ph I)$.  The shifted coronals are
\[
\begin{gathered}
\scoronal(W_{t,j})=
\begin{cases}
\dfrac{t}{t+\sig},&j=0,\\[3pt]
2+(t-3)\sig,&1\le j\le t,
\end{cases}\\[4pt]
\scoronal(C_5)=2-\sig,\quad
\scoronal(K_1\join P_3)=\scoronal(K_1\join P_4)=2,\\[2pt]
\scoronal(K_1\join C_5)=\scoronal(\overline{W_{3,3}})=2,\quad
\scoronal(K_1\join(K_2\udot K_2))=4+\sig,\quad
\scoronal(K_1\join(K_2\udot K_3))=12+6\sig.
\end{gathered}
\]

\end{lemma}
\begin{proof}
For $j=0$, $W_{t,0}=K_t$ and $A(K_t)+\ph I=J_t+\sig I_t\succ0$.
Since $(A(K_t)+\ph I)(t+\sig)^{-1}\bbone_t=\bbone_t$,
we have $\scoronal(K_t)=t/(t+\sig)$.

Suppose $1\le j\le t$.  After a simultaneous permutation of the rows
and columns of $A(W_{t,j})+\ph I$, with the clique vertices first, put
$B=\begin{pmatrix}I_j\\0\end{pmatrix}$.  Then
\[
A(W_{t,j})+\ph I=
\begin{pmatrix}
J_t+\sig I_t&B\\
B^{\mathsf T}&\ph I_j
\end{pmatrix}.
\]
Since $\ph I_j\succ0$, Lemma~\ref{lem:schur-congruence} gives
\[
\begin{aligned}
\operatorname{In}(A(W_{t,j})+\ph I)
&=\operatorname{In}(\ph I_j)
  +\operatorname{In}\bigl(J_t+\sig I_t-\sig BB^{\mathsf T}\bigr)\\
&=\operatorname{In}(\ph I_j)
  +\operatorname{In}\bigl(J_t+\sig\operatorname{diag}(0_j,I_{t-j})\bigr).
\end{aligned}
\]
Both matrices on the right are positive semidefinite, so
$A(W_{t,j})+\ph I\succeq0$.  

To verify the range condition, seek a solution of the form
\[
x=(a\bbone_j^{\mathsf T},\bbone_{t-j}^{\mathsf T},b\bbone_j^{\mathsf T})^{\mathsf T}.
\]
Multiplication by the displayed block matrix gives
\[
(A(W_{t,j})+\ph I)x=
\begin{pmatrix}
(ja+\sig a+t-j+b)\bbone_j\\
(ja+t-j+\sig)\bbone_{t-j}\\
(a+\ph b)\bbone_j
\end{pmatrix}.
\]
Choose $a=(\sig^2-t+j)/j$ and $b=\sig(1-a)$.
Then $ja+t-j=\sig^2$ and $a+\ph b=1$.
Using $\sig^2+\sig=1$, the first two block rows also equal the
corresponding all-ones vectors.  Thus
$(A(W_{t,j})+\ph I)x=\bbone$.
When $j=t$, the second block is empty and the same calculation applies.

Moreover,
$\scoronal(W_{t,j})=\bbone^{\mathsf T}x
=ja+(t-j)+jb
=2+(t-3)\sig$.

For the seven exceptional graphs, exact calculation verifies positive
semidefiniteness, the range condition, and the stated shifted-coronal
values.  See the \hyperref[sec:data-availability]{data availability statement}
for verification data.
\end{proof}

\begin{samepage}
\begin{definition}\label{def:admissible-class}
Call a graph $H$ \emph{admissible} if one of the following conditions
holds.  In \textup{(ii)}--\textup{(iv)}, the symbols $L$ and $R$
stand for disjoint unions of graphs in $\mathcal H_+$, possibly
empty, and $r=\scoronal(L)$, $s=\scoronal(R)$.
\begin{enumerate}[label=\textup{(\roman*)},leftmargin=*]
\item $H$ is a disjoint union of graphs in $\mathcal H_+$.
\item $H\cong(K_1\join L)\udot R$, where $R\ne\varnothing$, $r>\ph$,
and $s\le\sig+\frac{\sig^2}{r-\ph}$.
\item $H\cong\Sp(L,t)\udot R$ for some integer $t\ge1$, where
$R\ne\varnothing$ if $t=1$, $r>1$, and
$s\le\frac{\sig^2}{r-1}-\sig(t-2)$.
\item $H\cong\Spstar(L,t)\udot R$ for some integer $t\ge1$, where
$r+\sig(t-2)>1$ and
$s\le\frac{\sig^2}{r+\sig(t-2)-1}$.
\end{enumerate}
\end{definition}
\end{samepage}

\begin{theorem}\label{thm:global-clique-free}
Let $G$ be a connected graph with at least two vertices, and put
$H=\overline G$.  Then $\lambda_2(G)\le\sig$ if and only if $H$ is
admissible in the sense of Definition~\ref{def:admissible-class}.
\end{theorem}

The disconnected case is determined by this theorem as well.
If a graph $F$ has two components containing edges, it contains an
induced $2K_2$, so $\lambda_2(F)\ge\lambda_2(2K_2)=1>\sig$ by
Lemma~\ref{lem:induced-spectral-bounds}.  Otherwise, $F$ is edgeless
or has the form $C\udot kK_1$, where $C$ is connected of order at
least two and $k\ge0$.  In the latter case,
$\lambda_2(F)\le\sig$ if and only if $\lambda_2(C)\le\sig$,
since adding isolated vertices adds only zero eigenvalues and $\sig>0$.
Thus all graphs of order at least two satisfying the bound are
edgeless graphs or are obtained from the connected graphs in
Theorem~\ref{thm:global-clique-free} by adding isolated vertices.

In \textup{(ii)}--\textup{(iv)} of Definition~\ref{def:admissible-class},
the graphs forming $L$ and $R$ belong to the explicit family $\mathcal H_+$.
Their shifted coronals are given by Lemma~\ref{lem:positive-kernel-coronal},
and additivity determines $r$ and $s$.  Thus the numerical inequalities
can be checked from the component parameters without computing the
adjacency eigenvalues of $G$.  The proof of Theorem~\ref{thm:global-clique-free} is given in
Section~\ref{sec:global}.

\begin{example}
Let $L=4K_2$, $H_1=(K_1\join L)\udot K_4$, $H_2=(K_1\join L)\udot K_5$,
and $G_i=\overline{H_i}$ for $i=1,2$.
Both complements have the cone form in
Definition~\ref{def:admissible-class}\textup{(ii)}, with the same graph $L$
and with $R=K_4$ or $K_5$, respectively.  Their common value of $r$ is
\[
r=\scoronal(L)=4\scoronal(K_2)
=\frac{8}{2+\sig}=12-4\sqrt5>\ph .
\]
However, only $H_1$ satisfies the inequality in \textup{(ii)}:
\[
\scoronal(K_4)=\frac{4}{4+\sig}
<\sig+\frac{\sig^2}{r-\ph}
=\frac{33\sqrt5-19}{62}
<\frac{5}{5+\sig}=\scoronal(K_5).
\]
Consequently, $\lambda_2(G_1)\le\sig$, whereas
$\lambda_2(G_2)>\sig$.
\end{example}

\begin{remark}
The classes in Definition~\ref{def:admissible-class}\textup{(ii)}--\textup{(iv)}
need not be disjoint.  For $t=2$, thin and thick spiders are isomorphic,
and their numerical conditions coincide.  Also,
$\Spstar(L,1)\cong(K_1\join L)\udot K_1$, so the graphs in
\textup{(iv)} with $t=1$ are already included in \textup{(ii)}.
\end{remark}

\section{The family \texorpdfstring{$\mathcal H_+$}{H+} and forbidden induced subgraphs}\label{sec:positive-seeds}
We first prove that $\mathcal H_+$ consists exactly of the connected
graphs with least eigenvalue at least $-\ph$.
We then determine the minimal forbidden induced subgraphs for this
bound whose complements satisfy $\lambda_2\le\sig$.  These graphs
will be used to obtain the decompositions in
Sections~\ref{sec:a1} and~\ref{sec:clawfree}.

A graph $S$ is a \emph{minimal forbidden induced subgraph}
for the bound $\lambda_{\min}\ge-\ph$ if
$\lambda_{\min}(S)<-\ph$ and
$\lambda_{\min}(S-v)\ge-\ph$ for every $v\in V(S)$.
By interlacing, the vertex-deletion condition is equivalent to the
same bound holding for every nonempty proper induced subgraph of $S$.
Let $\mathcal F$ consist of those minimal forbidden induced subgraphs
$S$ that also satisfy $\lambda_2(\overline S)\le\sig$.

The following lemma is obtained by checking all nonempty connected
graphs of order at most seven.
See the \hyperref[sec:data-availability]{data availability statement}
for verification data.

\begin{lemma}\label{lem:positive-small-base}
Let $H$ be a connected graph with $\lambda_{\min}(H)\ge-\varphi$.
If $|V(H)|\le6$, then $H\in\mathcal H_+$.
If $|V(H)|=7$, then $H\in\mathcal W$.
\end{lemma}

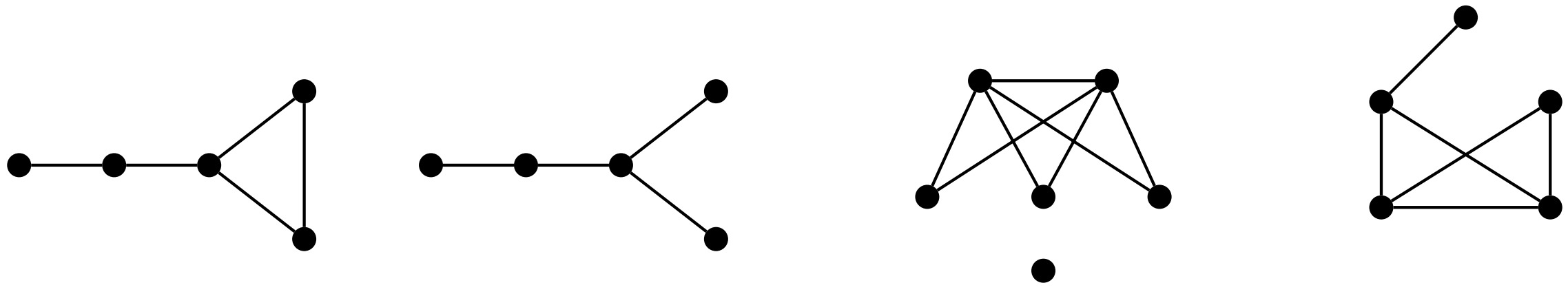
\begin{figure}[!htbp]
\centering
\resizebox{0.9\textwidth}{!}{%
\begin{tikzpicture}[
  line width=0.8pt,
  v/.style={circle,draw=black,fill=black,inner sep=0pt,minimum size=2mm},
  lab/.style={font=\small}
]
% ---------- F_1 ----------
\begin{scope}[shift={(0,0)}]
\node[v] (a) at (0,0) {};
\node[v] (b) at (0.9,0.7) {};
\node[v] (c) at (0.9,-0.7) {};
\node[v] (d) at (-0.9,0) {};
\node[v] (e) at (-1.8,0) {};
\draw (a)--(b) (a)--(c) (b)--(c) (a)--(d) (d)--(e);
\node[lab] at (0,-1.5) {$F_1$};
\end{scope}
% ---------- F_2 ----------
\begin{scope}[shift={(3.9,0)}]
\node[v] (a2) at (0,0) {};
\node[v] (b2) at (0.9,0.7) {};
\node[v] (c2) at (0.9,-0.7) {};
\node[v] (d2) at (-0.9,0) {};
\node[v] (e2) at (-1.8,0) {};
\draw (a2)--(b2) (a2)--(c2) (a2)--(d2) (d2)--(e2);
\node[lab] at (0,-1.5) {$F_2$};
\end{scope}
% ---------- F_3 ----------
\begin{scope}[shift={(7.9,0)}]
\node[v] (u) at (-0.6,0.8) {};
\node[v] (vv) at (0.6,0.8) {};
\node[v] (p) at (-1.1,-0.3) {};
\node[v] (qq) at (0,-0.3) {};
\node[v] (r) at (1.1,-0.3) {};
\node[v] (s) at (0,-1.0) {};
\draw (u)--(vv);
\draw (u)--(p) (u)--(qq) (u)--(r) (vv)--(p) (vv)--(qq) (vv)--(r);
\node[lab] at (0,-1.5) {$F_3$};
\end{scope}
% ---------- F_4 ----------
\begin{scope}[shift={(11.9,0)}]
\node[v] (uu) at (-0.8,0.6) {};
\node[v] (vv4) at (0.8,0.6) {};
\node[v] (w) at (-0.8,-0.4) {};
\node[v] (tt) at (0.8,-0.4) {};
\node[v] (pend) at (0,1.4) {};
\draw (uu)--(w) (uu)--(tt) (vv4)--(w) (vv4)--(tt) (w)--(tt) (pend)--(uu);
\node[lab] at (0,-1.5) {$F_4$};
\end{scope}
\end{tikzpicture}}
\caption{The forbidden induced subgraphs $F_1,\ldots,F_4$.}
\label{fig:f-graphs}
\end{figure}

Put $S_1=K_{1,3}$, $S_2=K_2\join(K_1\udot K_2)$,
$S_3=K_3\join(K_1\udot K_1)$, $S_4=K_1\join(K_2\udot K_4)$, and
$S_5=K_1\join(K_3\udot K_3)$.
Let $F_1$ be obtained by identifying one endpoint of a path of length
$2$ with a vertex of a triangle, and let $F_2$ be obtained by replacing
one edge of a claw with a path of length $2$ whose internal vertex is new.
Put $F_3=K_1\udot(K_2\join3K_1)$.  Let $F_4$ be obtained from $K_4$
by deleting one edge and attaching a leaf to a degree-two vertex of
the resulting graph; see Figure~\ref{fig:f-graphs}.

\FloatBarrier

The following lemma collects the spectral information on these graphs.
Part \textup{(i)} will be used in the least-eigenvalue classification,
and part \textup{(ii)} supplies the forbidden induced subgraphs for the
structural arguments.

\begin{lemma}\label{lem:fixed-forbidden-subgraphs}
\begin{enumerate}[label=\textup{(\roman*)},leftmargin=*]
\item $\lambda_{\min}(F)<-\ph$ for every
$F\in\{S_1,\ldots,S_5,C_4,F_1\}$.
\item If $G$ is connected, $H=\overline G$, and $\lambda_2(G)\le\sig$,
then $H$ contains none of $C_4,F_1,F_2,F_3,F_4$,
$S_i\udot S_j$, or $\Spstar(S_i,1)$ as an induced subgraph,
for $1\le i,j\le5$.
\end{enumerate}
\end{lemma}
\begin{proof}
Exact calculation gives
$n_-(A(F)+\ph I)\ge1$ for the graphs in part (i), and
$n_-(A(F)+\ph I-J)\ge2$ for every graph $F$ listed in part (ii).
Since $A(F)+\ph I-J=\sig I-A(\overline F)$, the latter condition
is equivalent to $\lambda_2(\overline F)>\sig$.
If such an $F$ occurred as an induced subgraph of $H$, then
$\lambda_2(\overline F)\le\lambda_2(\overline H)=\lambda_2(G)\le\sig$
by Lemma~\ref{lem:induced-spectral-bounds}, a contradiction.
See the \hyperref[sec:data-availability]{data availability statement}
for verification data.
\end{proof}

The finite checks above provide the base cases for the classification
of $\mathcal H_+$.  To extend it to arbitrary order, we show that a
connected graph obtained by adding a vertex to a graph in $\mathcal W$
of order at least seven either belongs to $\mathcal W$ or contains an
induced subgraph of order at most seven whose least eigenvalue is
below $-\ph$.

\begin{lemma}\label{lem:w-extension-dichotomy}
Suppose $H$ is connected, $H-x\cong W_{m,j}$, $0\le j\le m$, and
$|V(H)|\ge8$.  Then either $H\in\mathcal W$, or $H$ contains one of
$S_1,S_2,S_3,S_4,C_4,F_1$ as an induced subgraph.
\end{lemma}
\begin{proof}
Write $K=\{c_1,\ldots,c_m\}$ for the clique of $W_{m,j}$ and
$\ell_i$ for the leaf adjacent to $c_i$ ($1\le i\le j$).
Since $m+j\ge7$ and $j\le m$, we have $m\ge4$.
Suppose that $H$ contains none of $S_1,S_2,S_3,S_4,C_4,F_1$
as an induced subgraph.  We show that $H\in\mathcal W$.

Put $k=|N_H(x)\cap K|$.
If $2\le k\le m-2$, choose two neighbours and two non-neighbours
of $x$ in $K$; together with $x$, they induce $S_2$.
If $k=m-1$, choose three neighbours and the unique non-neighbour
of $x$ in $K$ to obtain $S_3$.
If $k=0$, connectedness gives $x\ell_i\in E(H)$ for some $i$;
then $x,\ell_i,c_i$ and two further vertices of $K$ induce $F_1$.
Hence $k=1$ or $k=m$.

If $k=m$, then $x$ has no neighbour among the leaves: otherwise
$x,\ell_i,c_i$ and two further vertices of $K$ induce $S_2$.
Consequently, $H\cong W_{m+1,j}$.

Suppose $k=1$, and let $c_h$ be the unique neighbour of $x$ in $K$.
If $h>j$, then $x$ has no neighbour among the leaves, since any edge
$x\ell_i$ would give an induced $C_4$ on $x,c_h,c_i,\ell_i$.
Thus $H\cong W_{m,j+1}$.

It remains to exclude $h\le j$.  Since
$\{c_h,x,\ell_h,c_i\}$ would induce $S_1$ for any $i\ne h$
if $x\ell_h\notin E(H)$, we must have $x\ell_h\in E(H)$.
Moreover, $x$ is nonadjacent to every $\ell_i$ with $i\ne h$,
since otherwise $x,c_h,c_i,\ell_i$ would induce $C_4$.
If $j\ge2$, choose $i\ne h$ with $1\le i\le j$; then
$x,\ell_h,c_h,c_i,\ell_i$ induce $F_1$.
If $j=1$, then $m+2=|V(H)|\ge8$, and the triangle
$\{x,\ell_h,c_h\}$ together with four vertices of $K\setminus\{c_h\}$
induces $S_4$.  Both cases contradict our assumption.
\end{proof}

\begin{proposition}
\label{prop:positive-complete}
A connected graph $H$ satisfies $\lambda_{\min}(H)\ge-\ph$
if and only if $H\in\mathcal H_+$.
\end{proposition}
\begin{proof}
By Lemma~\ref{lem:positive-kernel-coronal}, every graph in
$\mathcal H_+$ satisfies $\lambda_{\min}(H)\ge-\ph$.  Conversely,
Lemma~\ref{lem:positive-small-base} settles the case $|V(H)|\le7$.
Now suppose $|V(H)|\ge8$, and let $x$ be a leaf of a spanning tree of
$H$.  Then $H-x$ is connected and, by
Lemma~\ref{lem:induced-spectral-bounds},
$\lambda_{\min}(H-x)\ge-\ph$.  By induction on $|V(H)|$
(the case $|V(H-x)|=7$ is given by
Lemma~\ref{lem:positive-small-base}), $H-x\in\mathcal H_+$.
The seven exceptional graphs in $\mathcal H_+$ have order at most six,
so $H-x\in\mathcal W$.  Lemma~\ref{lem:w-extension-dichotomy} then
implies $H\in\mathcal W$; otherwise $H$ would contain one of
$S_1,S_2,S_3,S_4,C_4,F_1$ as an induced subgraph, contradicting
$\lambda_{\min}(H)\ge-\ph$ by
Lemma~\ref{lem:fixed-forbidden-subgraphs}\textup{(i)}.
\end{proof}

Proposition~\ref{prop:positive-complete} settles the case
$A(H)+\ph I\succeq0$.  We now consider graphs with
$\lambda_{\min}(H)<-\ph$ whose complements satisfy
$\lambda_2\le\sig$, and determine the family $\mathcal F$ introduced
at the beginning of the section.

\begin{lemma}\label{lem:seed-order}
Every member of $\mathcal F$ is connected and has at most seven
vertices.
\end{lemma}
\begin{proof}
Let $S\in\mathcal F$.  It is connected: otherwise a proper
component $C$ would have $\lambda_{\min}(C)<-\ph$, and for any
$v\notin C$ the induced subgraph $S-v$ would contain $C$, so
$\lambda_{\min}(S-v)<-\ph$, contradicting the definition of $\mathcal F$.
Suppose $|V(S)|\ge8$ and delete a leaf $x$ of a spanning tree of $S$.
Then $S-x$ is connected and, by minimality,
$\lambda_{\min}(S-x)\ge-\ph$.  Proposition~\ref{prop:positive-complete}
therefore gives $S-x\in\mathcal W$, since $|V(S-x)|\ge7$
and the seven exceptional graphs have order at most six.
Lemma~\ref{lem:w-extension-dichotomy} implies either
$S\in\mathcal W$ or that $S$ contains one of the listed induced
subgraphs $F$, each of order at most seven and with least eigenvalue
below $-\ph$ by Lemma~\ref{lem:fixed-forbidden-subgraphs}\textup{(i)}.
The first alternative contradicts the definition of $\mathcal F$, since
every member of $\mathcal W$ has least eigenvalue at least $-\ph$
by Lemma~\ref{lem:positive-kernel-coronal}.
In the second, $F$ is a proper induced subgraph of $S$, since
$|V(S)|\ge8$; choosing $v\in V(S)\setminus V(F)$ gives
$\lambda_{\min}(S-v)\le\lambda_{\min}(F)<-\ph$, again contradicting
the definition of $\mathcal F$.
\end{proof}

\begin{proposition}\label{prop:seed-classification}
Up to isomorphism, $\mathcal F=\{S_1,\ldots,S_5\}$.
\end{proposition}
\begin{proof}
By Lemma~\ref{lem:seed-order}, it suffices to consider connected graphs
of orders $4,\ldots,7$, since those of order at most three have least
eigenvalue at least $-\ph$.
Enumerating these graphs and imposing the three defining conditions
of $\mathcal F$ leaves precisely
$S_1,\ldots,S_5$, each once up to isomorphism.  The order bound in
Lemma~\ref{lem:seed-order} ensures that this finite enumeration is
exhaustive for $\mathcal F$.
See the \hyperref[sec:data-availability]{data availability statement}
for verification data.
\end{proof}

\begin{corollary}\label{cor:seed-detection}
If $\lambda_{\min}(H)<-\ph$ and $\lambda_2(\overline H)\le\sig$,
then $H$ contains an induced member of $\mathcal F$.
\end{corollary}
\begin{proof}
Choose an induced subgraph $S$ of $H$, minimal subject to
$\lambda_{\min}(S)<-\ph$.  By Lemma~\ref{lem:induced-spectral-bounds},
$\lambda_2(\overline S)\le \lambda_2(\overline H)\le \sig$, so $S\in\mathcal F$.
Proposition~\ref{prop:seed-classification} therefore gives $S\cong S_i$
for some $i\in\{1,\ldots,5\}$.
\end{proof}

The next lemma will be applied after a cone or spider decomposition
has been found.  It ensures that $L$ and $R$ satisfy
$A(L)+\ph I\succeq0$ and $A(R)+\ph I\succeq0$; by
Proposition~\ref{prop:positive-complete}, each of their connected
components therefore belongs to $\mathcal H_+$.

\begin{lemma}
\label{lem:constructor-positive-parts}
Let $G$ be a connected graph of order at least $2$, put $H=\overline G$,
and suppose $\lambda_2(G)\le\sig$.  Suppose $H=H_0\udot R$, where $H_0$ contains
an induced $S_j$ and has one of the forms
$K_1\join L$, $\Sp(L,t)$, or $\Spstar(L,t)$,
with $t\ge1$ in the latter two cases.  Then
\[
A(L)+\ph I\PSD,\qquad A(R)+\ph I\PSD.
\]
\end{lemma}
\begin{proof}
If $A(R)+\ph I\not\PSD$, then $\lambda_{\min}(R) < -\ph$ and
Lemma~\ref{lem:induced-spectral-bounds} gives
$\lambda_2(\overline R)\le \lambda_2(G)\le \sig$.  Corollary~\ref{cor:seed-detection}
produces an induced $S_i$ in $R$, so $H$ contains
$S_i\udot S_j$, contrary to
Lemma~\ref{lem:fixed-forbidden-subgraphs}\textup{(ii)}.

If $A(L)+\ph I\not\PSD$, then $\lambda_{\min}(L)< -\ph$ and
Lemma~\ref{lem:induced-spectral-bounds} gives
$\lambda_2(\overline L)\le\sig$.  Corollary~\ref{cor:seed-detection}
produces an induced $S_i$ in $L$.  If $H_0=\Spstar(L,t)$,
then $H_0$ contains an induced $\Spstar(S_i,1)$.
If $H_0=\Sp(L,t)$ with $t\ge2$, then $H_0$ contains
an induced $\Sp(L,2)\cong \Spstar(L,2)$,
and hence an induced $\Spstar(S_i,1)$.
If $H_0=K_1\join L$ or $H_0=\Sp(L,1)\cong K_1\join(L\udot K_1)$,
then $R\ne\varnothing$, since otherwise the apex would be
isolated in $G$, contradicting connectedness of $G$.  The apex and any vertex of $R$, together with
$S_i$, again induce $\Spstar(S_i,1)$.
Each case contradicts Lemma~\ref{lem:fixed-forbidden-subgraphs}\textup{(ii)},
proving $A(L)+\ph I\PSD$.
\end{proof}

\section{Structure of complements containing a claw}\label{sec:a1}
Throughout this section, $G$ is connected, $H=\overline G$, and
$\lambda_2(G)\le\sig$; we assume that $H$ contains an induced claw.
We show that the component of $H$ containing an induced claw is
a cone or a spider.
First, the forbidden induced subgraphs restrict how neighbours of its
centre can be adjacent to vertices outside its closed neighbourhood.
We then use these adjacency restrictions to obtain the partition in
Definition~\ref{def:spider-families}.

\begin{lemma}\label{lem:claw-dom-new}
Let $S'=H[\{c',a_1,a_2,a_3\}]$ be any induced claw of $H$ with centre $c'$.  Then
\[
N_H[a_i]\subseteq N_H[c'],\qquad i=1,2,3.
\]
Equivalently, if $x\notin V(S')$ and $xc'\notin E(H)$, then
$xa_i\notin E(H)$ for $i=1,2,3$.
\end{lemma}
\begin{proof}
Suppose $x\notin V(S')$, $xc'\notin E(H)$, and
$xa_i\in E(H)$ for some $i$.  If
$xa_i,xa_j\in E(H)$ for distinct $i,j$, then
$H[\{c',x,a_i,a_j\}]\cong C_4$.  If $x$ is adjacent to exactly one of
$a_1,a_2,a_3$, then $H[V(S')\cup\{x\}]\cong F_2$.  Both possibilities contradict
Lemma~\ref{lem:fixed-forbidden-subgraphs}.
\end{proof}

From now on, fix an induced claw
$S=H[\{c,\ell_1,\ell_2,\ell_3\}]$ with centre $c$.
By Lemma~\ref{lem:claw-dom-new}, no $\ell_i$ has a neighbour outside
$N_H[c]$.  Since $G$ is connected, $V(H)\setminus N_H[c]\ne\varnothing$;
otherwise $c$ would be isolated in $G$.

Define
\[
U=\{x\in N_H(c):N_H(x)\not\subseteq N_H[c]\},
\qquad D=N_H(U)\setminus N_H[c],
\]
and put
\[
O=V(H)\setminus(N_H[c]\cup D).
\]
Thus $D$ is the set of vertices at distance two from $c$ in $H$,
while $U$ consists of the vertices adjacent to $c$ and to at least
one vertex of $D$.
None of $\ell_1,\ell_2,\ell_3$ belongs to $U$, by
Lemma~\ref{lem:claw-dom-new}.
In particular, $D\cup O=V(H)\setminus N_H[c]$, and $O$ is
anticomplete to $N_H(c)$.  For $0\le i\le3$, let $U_i$ be the
set of vertices of $U$ adjacent to exactly $i$ vertices of
$\{\ell_1,\ell_2,\ell_3\}$.  For each $x\in U_2$,
let $\ell_x$ be the unique vertex of $\{\ell_1,\ell_2,\ell_3\}$
satisfying $x\ell_x\notin E(H)$.  Put $T_0=\{c\}\cup U_3$.
The sets $U$ and $D$ contain the endpoints of all edges with one endpoint
in $N_H[c]$ and the other outside it, whereas $O$ contains the remaining
vertices outside $N_H[c]$.  The next lemma controls these edges;
Lemma~\ref{lem:claw-rigidity} will construct a partition satisfying
the adjacency conditions in Definition~\ref{def:spider-families}.

\begin{lemma}\label{lem:claw-attachments}
\begin{enumerate}[label=\textup{(\roman*)}]
\item $U=U_2\cup U_3$, and $U_3$ is complete to $D\cup O$.
\item Every $x\in U_2$ has exactly one neighbour $p_x$ in $D$.
Moreover, if $w\in N_H(c)\setminus\{x,\ell_x\}$ is nonadjacent to $p_x$,
then $xw\in E(H)$ and $w\ell_x\notin E(H)$.
\item If $U_3\ne\varnothing$, then every $t\in T_0$ has a unique
non-neighbour $q_t$, and no vertex is nonadjacent to two members of
$T_0$.  Put $Q_0=\{q_t:t\in T_0\}$.  Then $|Q_0|=|T_0|$, and $Q_0$
is independent and anticomplete to $\{\ell_1,\ell_2,\ell_3\}$;
moreover, $q_t\in N_H(c)\setminus U$ for $t\in U_3$, and
$D=\{q_c\}$ and $O=\varnothing$.
\item $D$ is independent, and $O$ is anticomplete to
$N_H[c]\cup D$.
\item $U$ is a clique. If distinct $x,y\in U_2$ have a common
neighbour $p\in D$, then $\ell_x\ne\ell_y$ and $D=\{p\}$.
\end{enumerate}
\end{lemma}
\begin{proof}
For \textup{(i)}, let $x\in U_0\cup U_1$.  By the definition of $U$, there is a
$p\in D$ adjacent to $x$.  Choose distinct
$\ell_i,\ell_j\in\{\ell_1,\ell_2,\ell_3\}$ nonadjacent to $x$.
Then $H[\{p,x,c,\ell_i,\ell_j\}]\cong F_2$, a contradiction.
Thus $U=U_2\cup U_3$.  If $x\in U_3$ and $p\in D\cup O$
were nonadjacent, then $H[V(S)\cup\{p,x\}]\cong F_3$.
Hence $U_3$ is complete to $D\cup O$, proving \textup{(i)}.

For \textup{(ii)}, fix $x\in U_2$ and choose a neighbour
$p_x\in D$, which exists by the definition of $U$.  Write
$\{\ell_a,\ell_b\}=\{\ell_1,\ell_2,\ell_3\}\setminus\{\ell_x\}$.
By Lemma~\ref{lem:claw-dom-new}, $p_x$ is nonadjacent to
$\ell_a,\ell_b$, so $H[\{x,p_x,\ell_a,\ell_b\}]$
is a claw centred at $x$.
Suppose $x$ has another neighbour $p'\in D\setminus\{p_x\}$.
If $p'p_x\in E(H)$, then
$H[\{p_x,p',x,c,\ell_x\}]\cong F_1$; otherwise the same
vertices induce $F_2$.  Both possibilities are forbidden, so
$p_x$ is the unique neighbour of $x$ in $D$.

Now let $w\in N_H(c)\setminus\{x,\ell_x\}$ be nonadjacent to $p_x$.
For $w\in\{\ell_a,\ell_b\}$ the conclusion follows immediately.
Otherwise, suppose $wx\notin E(H)$.  Lemma~\ref{lem:claw-dom-new},
applied to the claw centred at $x$, gives
$w\ell_a,w\ell_b\notin E(H)$.
According as $w\ell_x\in E(H)$ or $w\ell_x\notin E(H)$,
$H[\{w,c,x,p_x,\ell_x\}]$ is isomorphic to $F_1$ or $F_2$,
contradicting Lemma~\ref{lem:fixed-forbidden-subgraphs}\textup{(ii)}.
Hence $wx\in E(H)$.  If $w\ell_x\in E(H)$, then
$H[\{x,c,w,\ell_x,p_x\}]\cong F_4$, another contradiction.
This proves \textup{(ii)}.

For \textup{(iii)}, assume $U_3\ne\varnothing$.
The set $T_0$ is a clique, since two nonadjacent
vertices of $T_0$ would induce $C_4$ together with
$\ell_1,\ell_2$.  Each $t\in T_0$ is the centre of an
induced claw on $\{t,\ell_1,\ell_2,\ell_3\}$.
By Lemma~\ref{lem:claw-dom-new}, every non-neighbour of $t$
is nonadjacent to the three leaves.  Moreover, no vertex can be
nonadjacent to two members of $T_0$, since it would induce $F_3$
with these members and the three leaves.

Every $t\in T_0$ has a non-neighbour because $G$ is connected.
Suppose that $t$ has two non-neighbours, say $a,b$.  Choose
$u\in T_0\setminus\{t\}$ and a non-neighbour $v$ of $u$.
Since each of $a,b,v$ is nonadjacent to a unique member of $T_0$
(namely $t$ for $a,b$, and $u$ for $v$), we have $au,bu,vt\in E(H)$.
If $ab\notin E(H)$, then
$H[\{u,t,a,b\}]$ is a claw centred at $u$.  Since
$vt\in E(H)$ but $vu\notin E(H)$, this contradicts
Lemma~\ref{lem:claw-dom-new}.  Thus $ab\in E(H)$.
If $va\in E(H)$, then $H[\{u,t,v,a\}]\cong C_4$, a contradiction.
Hence $va\notin E(H)$, and similarly $vb\notin E(H)$.
But then $H[\{v,t,u,a,b\}]\cong F_1$, another contradiction.
Hence each $t\in T_0$ has a unique non-neighbour $q_t$.

We have $|Q_0|=|T_0|$, since no vertex is nonadjacent to two
members of $T_0$.  The set $Q_0$ is independent: for distinct
$t,u\in T_0$, an edge $q_tq_u\in E(H)$ would give
$H[\{t,u,q_t,q_u\}]\cong C_4$.  It is anticomplete to the three
leaves by Lemma~\ref{lem:claw-dom-new}.

Since $q_c$ is the unique non-neighbour of $c$, we have
$D\cup O=\{q_c\}$.  As $U_3$ is complete to $D\cup O$ by
\textup{(i)}, it follows that $D=\{q_c\}$ and $O=\varnothing$.
For $t\in U_3$, the distinctness of $q_t$ and $q_c$ gives
$q_t\in N_H(c)$.  Since $Q_0$ is independent, $q_t$
is nonadjacent to $D=\{q_c\}$, so $q_t\notin U$.
This proves \textup{(iii)}.

For \textup{(iv)}, we show that every vertex of $D$ is
isolated in $H[D\cup O]$.  If $U_3\ne\varnothing$, then
$D\cup O=\{q_c\}$ by \textup{(iii)}, so this is immediate.
Otherwise, by \textup{(i)}, $U=U_2$.  Fix $p\in D$ and choose
$x\in U$ adjacent to $p$, as is possible by the definition of $D$.
Then $p=p_x$ by \textup{(ii)}.  Write
$\{\ell_a,\ell_b\}=\{\ell_1,\ell_2,\ell_3\}\setminus\{\ell_x\}$.
Since $p\notin N_H[c]$, Lemma~\ref{lem:claw-dom-new} shows that
$p$ is nonadjacent to $\ell_a,\ell_b$.  Thus
$H[\{x,p,\ell_a,\ell_b\}]$ is a claw centred at $x$.
Suppose $p$ has a neighbour $p'\in(D\cup O)\setminus\{p\}$.
Applying Lemma~\ref{lem:claw-dom-new} to this claw gives
$N_H[p]\subseteq N_H[x]$, and hence $xp'\in E(H)$.
This contradicts the uniqueness of $p$
in \textup{(ii)} if $p'\in D$, and the definition of $O$ if
$p'\in O$.  Thus $D$ is independent and anticomplete to $O$.
Since $O$ is anticomplete to $N_H[c]$ by definition, it is
anticomplete to $N_H[c]\cup D$, proving \textup{(iv)}.

For \textup{(v)}, we first show that $U$ is a clique.
Take distinct $x,y\in U$.  If they have a common neighbour
$p\in D$ and $xy\notin E(H)$, then
$H[\{c,x,p,y\}]\cong C_4$, a contradiction.
If they have no common neighbour in $D$, then
$x,y\notin U_3$ by \textup{(i)}.  Thus $x,y\in U_2$, so by \textup{(ii)} their unique
neighbours in $D$ are $p_x,p_y$; these are distinct, since
$x,y$ share no neighbour in $D$.  Hence $y$ is nonadjacent to
$p_x$, and as $y\in N_H(c)\setminus\{x,\ell_x\}$, \textup{(ii)}
gives $xy\in E(H)$.

Now let distinct $x,y\in U_2$ have a common neighbour
$p\in D$.  If $\ell_x=\ell_y$, then
$H[\{p,x,y,c,\ell_x\}]\cong F_4$, a contradiction.
Thus $\ell_x\ne\ell_y$.  Suppose there is a vertex
$p'\in D\setminus\{p\}$.  By \textup{(iii)}, $U_3=\varnothing$,
so $U=U_2$.  Choose $z\in N_H(p')\cap U$, which exists by
the definition of $D$.  By \textup{(ii)}, $z$ is nonadjacent to $p$.
Applying \textup{(ii)} to $x$ and $y$, each with $w=z$, shows that
$z$ is nonadjacent to both $\ell_x$ and $\ell_y$.
Since these leaves are distinct, this contradicts $z\in U_2$.
Consequently $D=\{p\}$.
\end{proof}

We now use these adjacency restrictions to construct a partition of $V(H)$.
When $|D|\ge2$, the unique neighbours in $D$ will give the matching of
a thin spider.  When $|D|=1$, the missing edges to the leaves $\ell_x$
and the vertices $q_t$ will give the matching in the complement of a
thick spider.
With $q_t$ as in Lemma~\ref{lem:claw-attachments}\textup{(iii)}, put
\[
\begin{aligned}
T&=\{c\}\cup U,\\
Q&=D\cup\{\ell_x:x\in U_2\}\cup\{q_t:t\in U_3\},\\
L&=H[N_H(c)\setminus(U\cup Q)].
\end{aligned}
\]
Recall that $T_0=\{c\}\cup U_3\subseteq T$.
\begin{lemma}\label{lem:claw-rigidity}
Suppose $D\ne\varnothing$.  Then $T$ is a clique, $Q$ is independent,
$V(H)=V(L)\udot Q\udot T\udot O$.
Furthermore, $V(L)$ is complete to $T$ and anticomplete to $Q$,
and $O$ is anticomplete to $V(L)\udot Q\udot T$.
Moreover, $|Q|=|T|$, and the following alternatives hold:
\begin{enumerate}[label=\textup{(\roman*)}]
\item If $|D|\ge2$, then $H[T,Q]\cong|T|K_2$.
\item If $|D|=1$, then $\overline H[T,Q]\cong|T|K_2$.
\end{enumerate}
\end{lemma}
\begin{proof}
We first verify that $T$ is a clique, $Q$ is independent, and $V(L)$
is complete to $T$ and anticomplete to $Q$.
We then determine the edges between $T$ and $Q$.
By
Lemma~\ref{lem:claw-attachments}\textup{(i),(iii)}, $N_H(c)$ is the disjoint union of $U$,
$\{\ell_x:x\in U_2\}$, $\{q_t:t\in U_3\}$, and $V(L)$, while
$D\cup O=V(H)\setminus N_H[c]$.  Hence the sets $T,Q,V(L),O$
partition $V(H)$.  By Lemma~\ref{lem:claw-attachments}\textup{(v)},
$U$ is a clique, so $T=\{c\}\cup U$ is a clique.

The set $D$ is independent by
Lemma~\ref{lem:claw-attachments}\textup{(iv)}, and
$\ell_1,\ell_2,\ell_3$ are independent by the definition of $S$.
If $U_3=\varnothing$, then $Q=D\cup\{\ell_x:x\in U_2\}$
is independent, since Lemma~\ref{lem:claw-dom-new}
shows that $D$ is anticomplete to the three leaves.
If $U_3\ne\varnothing$, then Lemma~\ref{lem:claw-attachments}\textup{(iii)}
gives $D=\{q_c\}$ and shows that the vertices $q_t$ for
$t\in\{c\}\cup U_3$ are independent and anticomplete to the three leaves.
Hence $Q$ is independent in either case.

Let $w\in V(L)$.  Since $w\in N_H(c)\setminus U$,
it has no neighbours in $D\cup O$ by the definition of $U$.
For each $x\in U_2$, we have $w\ne\ell_x$ since $w\notin Q$,
so Lemma~\ref{lem:claw-attachments}\textup{(ii)} gives
$wx\in E(H)$ and $w\ell_x\notin E(H)$.
For $t\in U_3$, the unique non-neighbour $q_t$ belongs to $Q$,
so $wt\in E(H)$.  Moreover, $tq_c\in E(H)$ by
Lemma~\ref{lem:claw-attachments}\textup{(i)}, while
$wq_c\notin E(H)$.
If $wq_t\in E(H)$, then
$H[\{c,t,q_c,q_t,w\}]\cong F_4$, a contradiction.
Thus $wq_t\notin E(H)$.
Thus $V(L)$ is complete to $T$ and anticomplete to $Q$.
By Lemma~\ref{lem:claw-attachments}\textup{(iv)}, $O$ is
anticomplete to $N_H[c]\cup D$, hence to $V(L)\udot Q\udot T$.

It remains to determine the structure of $H[T,Q]$.

\textup{(i)} Suppose $|D|\ge2$.  By
Lemma~\ref{lem:claw-attachments}\textup{(iii)}, we have
$U_3=\varnothing$, so $U=U_2$.  Each $x\in U$ has exactly
one neighbour $p_x$ in $D$ by
Lemma~\ref{lem:claw-attachments}\textup{(ii)}, while every
$p\in D$ has a neighbour in $U$ by the definition of $D$.
Since $|D|\ge2$, Lemma~\ref{lem:claw-attachments}\textup{(v)}
shows that no vertex of $D$ can have two neighbours in $U$.
Hence every vertex of $D$ has exactly one neighbour in $U$.

If $x,y\in U$ are distinct, then $p_x\ne p_y$ and
Lemma~\ref{lem:claw-attachments}\textup{(ii)}, applied to $x$
with $w=y$, gives $y\ell_x\notin E(H)$.  Thus
$\ell_x=\ell_y$, since $y\in U_2$.  Write $\ell$ for this
common leaf.  We have $Q=D\cup\{\ell\}$ and
$|Q|=|D|+1=|U|+1=|T|$.  The only edges between $T$ and
$Q$ are $c\ell$ and $xp_x$ for $x\in U$, so
$H[T,Q]\cong|T|K_2$.  Together with the adjacency properties proved
above, this gives $H[V(L)\cup T\cup Q]\cong\Sp(L,|T|)$.

\textup{(ii)} Suppose $D=\{p\}$.  By the definition of $U$,
every vertex of $U$ is adjacent to $p$, whereas
$cp\notin E(H)$.  By Lemma~\ref{lem:claw-attachments}\textup{(v)},
the vertices $\ell_x$ for $x\in U_2$ are pairwise distinct.  Hence each such $\ell_x$
is adjacent to every vertex of $T$ except $x$, since every
$y\in U_2\setminus\{x\}$ is nonadjacent only to $\ell_y$ among the leaves
and every $t\in U_3$ is adjacent to all three leaves.
For $t\in U_3$, Lemma~\ref{lem:claw-attachments}\textup{(iii)}
gives a distinct vertex $q_t$.  It is adjacent to every vertex
of $T$ except $t$: by \textup{(iii)} we have $q_t\in N_H(c)$, and
no vertex is nonadjacent to two members of $T_0$, so $q_t$ is
adjacent to $c$ and to every member of $U_3\setminus\{t\}$.
If $q_t$ were nonadjacent to some $x\in U_2$, then the vertices
$c,t,x,q_c,q_t$ would induce $F_4$: indeed $c,t,x,q_c$ induce
$K_4-e$ with the unique nonedge $cq_c$, and $q_t$ is adjacent
only to $c$ among these vertices; this contradicts
Lemma~\ref{lem:fixed-forbidden-subgraphs}\textup{(ii)}.
Thus $q_t$ is adjacent to every vertex of $T$ except $t$;
these vertices are distinct from $p$ and from the leaves $\ell_x$.
Consequently, $|Q|=1+|U_2|+|U_3|=|T|$, and the only
nonedges between $T$ and $Q$ are $cp$, $x\ell_x$ for
$x\in U_2$, and $tq_t$ for $t\in U_3$, so
$\overline H[T,Q]\cong|T|K_2$.  Together with the adjacency properties proved above,
this gives $H[V(L)\cup T\cup Q]\cong\Spstar(L,|T|)$.
\end{proof}

\begin{theorem}\label{thm:A1-complete}
Let $G$ be connected, $H=\overline G$, and $\lambda_2(G)\le\sig$.
Suppose that $H$ contains an induced claw.  Then there exist graphs
$L$ and $R$ with $A(L)+\ph I\PSD$ and $A(R)+\ph I\PSD$ such that
$H$ is isomorphic to one of
\[
(K_1\join L)\udot R,\qquad
\Sp(L,t)\udot R,\qquad \Spstar(L,t)\udot R,
\]
where $t\ge1$ in the spider cases.
\end{theorem}
\begin{proof}
We distinguish three cases according to $D$, and put $R=H[O]$.

If $D=\varnothing$, then $U=\varnothing$, and
$O=V(H)\setminus N_H[c]$ is anticomplete to $N_H(c)$.
Since $G$ is connected, $O\ne\varnothing$; otherwise $c$ would be
isolated in $G$.
Taking $L=H[N_H(c)]$ gives
$H=(K_1\join L)\udot R$.

If $|D|\ge2$, Lemma~\ref{lem:claw-rigidity}\textup{(i)}
gives $H=\Sp(L,|T|)\udot R$, where $L,T,R$ are as above.  If $|D|=1$, Lemma~\ref{lem:claw-rigidity}\textup{(ii)} gives
$H=\Spstar(L,|T|)\udot R$.

In each case, the cone $K_1\join L$, $\Sp(L,|T|)$,
or $\Spstar(L,|T|)$ contains the induced
claw $S$.  Lemma~\ref{lem:constructor-positive-parts}
therefore yields $A(L)+\ph I\PSD$ and $A(R)+\ph I\PSD$.
\end{proof}

\section{Structure of claw-free complements}\label{sec:clawfree}

Assume throughout this section that $G$ is connected, $H=\overline G$ is
claw-free, and $\lambda_2(G)\le\sig$.  For an induced copy $S$ of any of the graphs $S_2,\ldots,S_5$,
there is a partition $V(S)=X\cup Y\cup Z$ such that
$S=H[X]\join(H[Y]\udot H[Z])$,
where $X,Y,Z$ induce cliques and $Y$ is anticomplete to $Z$.
For $S=S_2,S_3,S_4,S_5$, the triples $(|X|,|Y|,|Z|)$ are
$(2,1,2)$, $(3,1,1)$, $(1,2,4)$, and $(1,3,3)$, respectively.  Thus each vertex of $X$ is universal in $S$.

For $v\notin V(S)$, write $N_S(v)=N_H(v)\cap V(S)$, and define
\[
O_S=\{v\notin V(S):N_S(v)=\varnothing\},\qquad
U_S=\{v\notin V(S):N_S(v)\ne\varnothing\}.
\]
Thus the vertices outside $S$ are partitioned according as they have
no neighbour in $S$ or at least one neighbour in $S$.

\begin{lemma}\label{lem:clawfree-attachment}
If $S\in\{S_2,S_4,S_5\}$, then $U_S$ is complete to $X$.
\end{lemma}
\begin{proof}
First let $S=S_2$, with $X=\{x_1,x_2\}$, $Y=\{y\}$, and
$Z=\{z_1,z_2\}$.  Suppose that $v\in U_S$ is not adjacent to $x_1$.
If $v$ has neighbours in both $Y$ and $Z$, then
$H[\{v,y,x_1,z_k\}]\cong C_4$ for some $z_k\in Z$,
contradicting Lemma~\ref{lem:fixed-forbidden-subgraphs}\textup{(ii)}.
Thus $v$ has neighbours in at most one of $Y,Z$.
If $vx_2\in E(H)$, then, since $H$ is claw-free, $v$ is adjacent
to $y$ or to both $z_1,z_2$; otherwise $\{x_2,v,y,z_k\}$ induces
a claw for some $z_k\in Z$.  Consequently, $N_H(v)\cap\{y,z_1,z_2\}\in\bigl\{\{y\},\{z_1,z_2\}\bigr\}$.
In these two cases, $H[\{v,x_1,y,z_1,z_2\}]$ is isomorphic
to $F_1$ or $F_4$, respectively, a contradiction.
If $vx_2\notin E(H)$, then by the definition of $U_S$ there
exists $u\in N_H(v)\cap(Y\cup Z)$.  Choose $u'$ from the other
of the two cliques $Y,Z$: if $u\in Y$, take $u'\in Z$, and conversely.
Since $v$ has neighbours in at most
one of $Y,Z$, we have $vu'\notin E(H)$, and
$H[\{u,u',x_1,x_2,v\}]\cong F_4$.
Both cases contradict Lemma~\ref{lem:fixed-forbidden-subgraphs}\textup{(ii)}.
Thus $v$ is adjacent to $x_1$, and by symmetry to $x_2$.

Next let $S=S_4$ or $S_5$, with $X=\{x\}$ and $|Y|,|Z|\ge2$.
Suppose that $v\in U_S$ is not adjacent to $x$.  By the definition
of $U_S$, the vertex $v$ has a neighbour in $Y\cup Z$.
If $v$ had neighbours $y\in Y$ and $z\in Z$, then $H[\{v,y,x,z\}]\cong C_4$, contradicting
Lemma~\ref{lem:fixed-forbidden-subgraphs}\textup{(ii)}.
Suppose its neighbours lie in $Y$.
If $v$ has exactly one neighbour $y_i\in Y$, then $v,y_i,x$
and two vertices of $Z$ induce $F_1$; if it has at least two neighbours
$y_i,y_j\in Y$, then $v,y_i,y_j,x$ and any vertex of $Z$ induce $F_4$.
The same argument applies with $Y,Z$ interchanged.  Both alternatives
contradict Lemma~\ref{lem:fixed-forbidden-subgraphs}\textup{(ii)}.
Hence $vx\in E(H)$.
\end{proof}

\begin{lemma}\label{lem:clawfree-separation-new}
If $S\in\{S_2,S_4,S_5\}$, then $O_S\ne\varnothing$ and
$O_S$ is anticomplete to $U_S$.
\end{lemma}
\begin{proof}
If $O_S=\varnothing$, Lemma~\ref{lem:clawfree-attachment}
makes every vertex of $X$ universal in $H$, hence isolated in $G$,
a contradiction to connectedness.
Let $o\in O_S$ and $v\in U_S$, and suppose for a contradiction
that $ov\in E(H)$.
The vertex $v$ is not complete to $V(S)$, since otherwise
$H[\{v,o,y,z\}]$ is a claw for $y\in Y$, $z\in Z$.

Suppose $S=S_2$.  By Lemma~\ref{lem:clawfree-attachment}, $v$ is
complete to $X$, so there is $u\in Y\cup Z$ with $vu\notin E(H)$.
Then $H[X\cup\{o,v,u\}]\cong F_4$, a contradiction.

For $S=S_4$ or $S_5$, write $X=\{x\}$.
Since $vx\in E(H)$ and $H$ is claw-free, $v$ is complete to at
least one of $Y,Z$: otherwise a non-neighbour from each set,
together with $x,v$, would induce a claw.
By symmetry, suppose that $v$ is complete to $Y$.
Since $v$ is not complete to $V(S)$ but is complete to $X\cup Y$,
there is $z_1\in Z$ with $vz_1\notin E(H)$.
Choose $z_2\in Z\setminus\{z_1\}$, which is possible because $|Z|\ge2$.
Then $H[\{x,z_1,z_2,v,o\}]\cong F_1$ or $F_4$ according as
$vz_2\notin E(H)$ or $vz_2\in E(H)$, contradicting
Lemma~\ref{lem:fixed-forbidden-subgraphs}\textup{(ii)}.
Hence $O_S$ is anticomplete to $U_S$.
\end{proof}

\begin{samepage}
\begin{lemma}\label{lem:clawfree-s3}
Suppose $S$ is an induced copy of $S_3$ in $H$, and $H$ has no induced $S_2$.
Then $O_S\ne\varnothing$ and $U_S=\varnothing$.
\end{lemma}
\begin{proof}
Write $X=\{x_1,x_2,x_3\}$, $Y=\{y\}$, and $Z=\{z\}$.
We first show that every $v\in U_S$ is complete to $V(S)$.
If $v$ is adjacent to neither $y$ nor $z$, then by the definition
of $U_S$ we may assume $vx_1\in E(H)$, and
$H[\{x_1,y,z,v\}]$ is a claw, a contradiction.
If $vy,vz\in E(H)$ but $vx_i\notin E(H)$ for some $i$, then
$H[\{v,y,x_i,z\}]\cong C_4$.
Finally, suppose $vy\in E(H)$ and $vz\notin E(H)$; the reverse
case is symmetric.  If $v$ has at most one neighbour in $X$,
choose two distinct non-neighbours $x_i,x_j\in X$ of $v$.
Then $H[\{v,y,z,x_i,x_j\}]\cong F_4$.
If $v$ has at least two neighbours $x_i,x_j\in X$, then
$H[\{v,y,z,x_i,x_j\}]\cong S_2$.
These possibilities are excluded by Lemma~\ref{lem:fixed-forbidden-subgraphs}\textup{(ii)}
and the assumption that $H$ contains no $S_2$.
Thus every $v\in U_S$ is complete to $V(S)$.

If $O_S=\varnothing$, every vertex of $X$ is universal in $H$ and
thus isolated in $G$, contrary to connectedness.  Fix $o\in O_S$.
Suppose $v\in U_S$.  If $ov\in E(H)$, then
$H[\{v,o,y,z\}]$ is a claw, contrary to the assumption that
$H$ is claw-free.  Otherwise,
$H[V(S)\cup\{o,v\}]\cong\Spstar(S,1)$, contradicting
Lemma~\ref{lem:fixed-forbidden-subgraphs}\textup{(ii)}.
Consequently $U_S=\varnothing$.
\end{proof}
\end{samepage}

\begin{theorem}\label{thm:clawfree-closure-new}
If $H$ contains an induced $S_i$ for some $i\in\{2,3,4,5\}$, then
there exist graphs $L$ and $R$, with $R\ne\varnothing$, such that
\[
H=(K_1\join L)\udot R,\qquad
A(L)+\ph I\PSD,\qquad A(R)+\ph I\PSD.
\]
\end{theorem}
\begin{proof}
If $H$ contains an induced $S_i$ for some $i\in\{2,4,5\}$,
choose such an $S$; otherwise, choose an induced $S_3$.
In the first case, Lemmas~\ref{lem:clawfree-attachment} and
\ref{lem:clawfree-separation-new} apply, whereas in the
second case, Lemma~\ref{lem:clawfree-s3} gives
$U_S=\varnothing$ and $O_S\ne\varnothing$.
Thus, in either case, $O_S\ne\varnothing$ and
$O_S$ is anticomplete to $V(S)\cup U_S$, with $U_S$
complete to $X$ whenever $U_S\ne\varnothing$.
Choose $c\in X$ and let
$L=H[(V(S)\cup U_S)\setminus\{c\}]$ and $R=H[O_S]$.
Then $H=(K_1\join L)\udot R$.  Since $K_1\join L$
contains an induced $S$, Lemma~\ref{lem:constructor-positive-parts}
yields $A(L)+\ph I\PSD$ and $A(R)+\ph I\PSD$.
Since $O_S\ne\varnothing$, we have $R\ne\varnothing$.
\end{proof}

\section{Proof of the main theorem}\label{sec:global}
We now derive the numerical conditions for the cone and spider
decompositions obtained in Sections~\ref{sec:a1} and~\ref{sec:clawfree}.
We first give a criterion for the rank-one perturbation $A(H)+\ph I-J$,
then compute the inertia and shifted coronal of each construction,
and use these results to prove Theorem~\ref{thm:global-clique-free}.

\begin{lemma}\label{lem:spectral-coronal-new}
Let $G$ be a connected graph of order at least $2$, and put $H=\overline G$.
If $A(H)+\ph I\succeq0$, then $\lambda_2(G)\le\sig$.
If $A(H)+\ph I\not\succeq0$ and
$\bbone\in\operatorname{range}(A(H)+\ph I)$, then
$\lambda_2(G)\le\sig$ if and only if
$n_-(A(H)+\ph I)=1$ and $\scoronal(H)\le1$.
\end{lemma}
\begin{proof}
Since $A(G)=J-I-A(H)$, Lemma~\ref{lem:weyl} gives, when
$A(H)+\ph I\succeq0$,
\[
\lambda_2(G)\le\lambda_2(J)+\lambda_1(-I-A(H))
=-1-\lambda_{\min}(H)\le\sig,
\]
where $\lambda_2(J)=0$, since $|V(G)|\ge2$.
For the second assertion, $\lambda_1(G)\ge1>\sig$ by
connectedness, and $A(H)+\ph I-J=\sig I-A(G)$.  Hence
$\lambda_2(G)\le\sig$ if and only if
$n_-(A(H)+\ph I-J)=1$.
Put $M=A(H)+\ph I$ and consider the matrix
\[
K=\begin{pmatrix}
M&\bbone\\
\bbone^{\mathsf T}&1
\end{pmatrix}.
\]
Taking the Schur complement of the scalar block $1$ gives
$n_-(K)=n_-(M-J)$.  On the other hand, the assumed range condition
allows us to take the generalized Schur complement of $M$, giving
\[
n_-(K)=n_-(M)+n_-(1-\bbone^{\mathsf T}M^\dagger\bbone)
=n_-(M)+n_-(1-\scoronal(H)).
\]
Consequently,
\[
n_-(M-J)=n_-(M)+n_-(1-\scoronal(H)).
\]
Since $M$ is not positive semidefinite, $n_-(M)\ge1$.
The left-hand side equals one precisely when $n_-(M)=1$ and
$1-\scoronal(H)\ge0$, proving the equivalence.
\end{proof}

\begin{samepage}
\begin{lemma}\label{lem:negative-coronal-calculus}
Let $L$ be a disjoint union of graphs in $\mathcal H_+$, possibly
empty, put $r=\scoronal(L)$, and let $t\ge1$ be an integer.
\begin{enumerate}[label=(\roman*)]
\item The negative inertia indices satisfy
\begin{align*}
n_-(A(K_1\join L)+\ph I)&\le1,
&&\text{with equality iff }r>\ph;\\
n_-(A(\Sp(L,t))+\ph I)&\le1,
&&\text{with equality iff }r>1;\\
n_-(A(\Spstar(L,t))+\ph I)&\le1,
&&\text{with equality iff }r+\sig(t-2)>1.
\end{align*}
\item For each $H_0\in\{K_1\join L,\Sp(L,t),\Spstar(L,t)\}$,
if the denominator in the corresponding formula below is nonzero,
then $\bbone\in\operatorname{range}(A(H_0)+\ph I)$, and
$\scoronal(H_0)$ is given by that formula:
\begin{align}
\scoronal(K_1\join L)
&=r+\frac{(1-r)^2}{\ph-r},
\label{eq:apex-coronal}\\
\scoronal(\Sp(L,t))
&=1+\sig(t-2)-\frac{\sig^2}{r-1},
\label{eq:thin-coronal}\\
\scoronal(\Spstar(L,t))
&=1-\frac{\sig^2}{r+\sig(t-2)-1}.
\label{eq:thick-coronal}
\end{align}
The strict inequalities in \textup{(i)} ensure that the corresponding
denominators in \textup{(ii)} are nonzero, so these formulas apply in
every case with one negative eigenvalue, even when the shifted
adjacency matrix is singular.
\end{enumerate}
\end{lemma}
\end{samepage}
\begin{proof}
Let $H_0$ be either $\Sp(L,t)$ or $\Spstar(L,t)$, and let
$V(H_0)=V(L)\cup T\cup Q$ be the partition in
Definition~\ref{def:spider-families}.  Denote by $B$ the block of
$A(H_0)$ with rows indexed by $T$ and columns by $Q$.
With the labelling described in Section~\ref{sec:matrix}, $B=I_t$
for the thin spider and $B=J_t-I_t$ for the thick spider.
Set $k=1$ in the first case and $k=t-1$ in the second.
In both cases,
$B\bbone_t=B^{\mathsf T}\bbone_t=k\bbone_t$ and
$BB^{\mathsf T}=I_t+(k-1)J_t$.

Put $n=|V(L)|$.  By Lemma~\ref{lem:positive-kernel-coronal},
$A(L)+\ph I\succeq0$ and
$\bbone_n\in\operatorname{range}(A(L)+\ph I)$.
Moreover, there exists $y$ satisfying
$(A(L)+\ph I)y=\bbone_n$ and $\bbone_n^{\mathsf T}y=r$.
With the rows and columns indexed by $V(L),T,Q$ in this order, we have
\[
A(H_0)+\ph I=
\begin{pmatrix}
A(L)+\ph I&J_{n,t}&0\\
J_{t,n}&J_t+\sig I_t&B\\
0&B^{\mathsf T}&\ph I_t
\end{pmatrix}.
\]
Since $J_{n,t}=\bbone_n\bbone_t^{\mathsf T}$, every column of
$\begin{pmatrix}J_{n,t}&0\end{pmatrix}$ belongs to
$\operatorname{range}(A(L)+\ph I)$.  Also,
\[
J_{t,n}(A(L)+\ph I)^\dagger J_{n,t}
=\bigl(\bbone_n^{\mathsf T}(A(L)+\ph I)^\dagger\bbone_n\bigr)J_t
=rJ_t.
\]
Thus Lemma~\ref{lem:schur-congruence} applies to the first block,
even if $A(L)+\ph I$ is singular, and gives
\[
\begin{aligned}
\operatorname{In}(A(H_0)+\ph I)
={}&\operatorname{In}(A(L)+\ph I)
+\operatorname{In}\!\begin{pmatrix}
(1-r)J_t+\sig I_t&B\\
B^{\mathsf T}&\ph I_t
\end{pmatrix}\\
={}&\operatorname{In}(A(L)+\ph I)
+\operatorname{In}(\ph I_t)+\operatorname{In}\bigl((1-r)J_t+\sig I_t
-\sig BB^{\mathsf T}\bigr)\\
={}&\operatorname{In}(A(L)+\ph I)
+\operatorname{In}(\ph I_t)
+\operatorname{In}\bigl([1-r-\sig(k-1)]J_t\bigr).
\end{aligned}
\]
For the second equality we take the Schur complement of the
positive definite block $\ph I_t$, using $\ph^{-1}=\sig$.
Since $J_t$ has exactly one positive eigenvalue,
$A(H_0)+\ph I$ is positive semidefinite precisely when
$1-r-\sig(k-1)\ge0$; otherwise it has exactly one negative
eigenvalue.  Substituting $k=1$ and $k=t-1$ proves \textup{(i)} for the
two spider families.

Suppose $1-r-\sig(k-1)\ne0$.  Set
\[
v=\frac{1-r-\sig k}{t[1-r-\sig(k-1)]}\bbone_t,\qquad
u=(1-\bbone_t^{\mathsf T}v)y,\qquad
w=\sig(\bbone_t-B^{\mathsf T}v),
\]
and put $x=(u^{\mathsf T},v^{\mathsf T},w^{\mathsf T})^{\mathsf T}$.

The first and third block rows give
$(A(L)+\ph I)u+J_{n,t}v=\bbone_n$ and
$B^{\mathsf T}v+\ph w=\bbone_t$.
For the second block row, using $\bbone_n^{\mathsf T}y=r$,
$B\bbone_t=k\bbone_t$ and
$BB^{\mathsf T}=I_t+(k-1)J_t$, we obtain
\[
\begin{aligned}
J_{t,n}u+(J_t+\sig I_t)v+B w
&=r(1-\bbone_t^{\mathsf T}v)\bbone_t
  +J_tv+\sig v+\sig B\bbone_t-\sig BB^{\mathsf T}v\\
&=(r+\sig k)\bbone_t
 +[1-r-\sig(k-1)]J_tv\\
&=\bbone_t.
\end{aligned}
\]
Thus $(A(H_0)+\ph I)x=\bbone$, so
$\bbone\in\operatorname{range}(A(H_0)+\ph I)$ and
\[
\begin{aligned}
\scoronal(H_0)
&=\bbone_n^{\mathsf T}u+\bbone_t^{\mathsf T}v
  +\bbone_t^{\mathsf T}w\\
&=r+\sig t+(1-r-\sig k)\bbone_t^{\mathsf T}v\\
&=1+\sig(t-k-1)+\frac{\sig^2}{1-r-\sig(k-1)}.
\end{aligned}
\]
Taking $k=1$ and $k=t-1$ yields
\eqref{eq:thin-coronal} and \eqref{eq:thick-coronal}.

Finally, consider the cone $K_1\join L$.
Since $\Spstar(L,1)\cong(K_1\join L)\udot K_1$, the extra isolated
vertex contributes only the positive eigenvalue $\ph$ to the shifted
adjacency matrix.  The thick-spider result with $t=1$ therefore gives
\[
A(K_1\join L)+\ph I\succeq0
\quad\Longleftrightarrow\quad
r-\sig\le1
\quad\Longleftrightarrow\quad
r\le\ph.
\]
For $r>\ph$, the shifted adjacency matrix of the cone has exactly
one negative eigenvalue.  This proves \textup{(i)} for cones.
If $r\ne\ph$, the result for $\Spstar(L,1)$ also shows that
$\bbone\in\operatorname{range}(A(K_1\join L)+\ph I)$.
Since $\scoronal(K_1)=\sig$, additivity gives
\[
\begin{aligned}
\scoronal(K_1\join L)
&=\scoronal(\Spstar(L,1))-\scoronal(K_1)\\
&=1-\sig+\frac{\sig^2}{\ph-r}
 =r+\frac{(1-r)^2}{\ph-r},
\end{aligned}
\]
which completes the proof.
\end{proof}

\begin{proof}[Proof of Theorem~\ref{thm:global-clique-free}]
We first prove necessity.  Suppose $\lambda_2(G)\le\sig$, and put
$H=\overline G$.
If $A(H)+\ph I\succeq0$, each connected component of $H$ satisfies
the same bound.  Proposition~\ref{prop:positive-complete}, applied to
each component, shows that $H$ satisfies
Definition~\ref{def:admissible-class}\textup{(i)}.

Now suppose $A(H)+\ph I$ is not positive semidefinite.
Corollary~\ref{cor:seed-detection} gives an induced $S_i$ in $H$.
If $H$ contains a claw, Theorem~\ref{thm:A1-complete} gives a decomposition
$H=H_0\udot R$, where $H_0$ is a cone or spider on $L$ and
$A(L)+\ph I\succeq0$, $A(R)+\ph I\succeq0$.
If $H$ is claw-free, then $i\in\{2,3,4,5\}$, and
Theorem~\ref{thm:clawfree-closure-new} gives the same conclusion with
$H_0=K_1\join L$ and $R\ne\varnothing$.
Thus, in either case, Proposition~\ref{prop:positive-complete}, applied
componentwise, shows that $L$ and $R$ are disjoint unions of members
of $\mathcal H_+$.
Put $r=\scoronal(L)$ and $s=\scoronal(R)$.
In the cone case, $R\ne\varnothing$, since otherwise the apex would
be isolated in $G$.  The same argument applies to the thin spider
with $t=1$, which is a cone over $L\udot K_1$.

Because $A(H)+\ph I$ is block diagonal with blocks
$A(H_0)+\ph I$ and $A(R)+\ph I$, and the latter is positive
semidefinite, $A(H_0)+\ph I$ is not positive semidefinite.
Lemma~\ref{lem:negative-coronal-calculus} now gives the corresponding
strict condition on $r$ in Definition~\ref{def:admissible-class}\textup{(ii)}--\textup{(iv)}.
It also gives $n_-(A(H_0)+\ph I)=1$ and
$\bbone\in\operatorname{range}(A(H_0)+\ph I)$.
For every component of $R$, Lemma~\ref{lem:positive-kernel-coronal}
gives positive semidefiniteness and the same range condition.
The block-diagonal structure therefore yields
\[
n_-(A(H)+\ph I)=1,\qquad
\bbone\in\operatorname{range}(A(H)+\ph I).
\]
By Lemma~\ref{lem:spectral-coronal-new} and additivity,
\[
\scoronal(H)=\scoronal(H_0)+s\le1.
\]
For a cone, formula~\eqref{eq:apex-coronal} rewrites this as
\[
s+1-\sig-\frac{\sig^2}{r-\ph}\le1
\quad\Longleftrightarrow\quad
s\le\sig+\frac{\sig^2}{r-\ph}.
\]
For thin and thick spiders, formulas~\eqref{eq:thin-coronal}
and~\eqref{eq:thick-coronal} give, respectively,
\[
s+1+\sig(t-2)-\frac{\sig^2}{r-1}\le1
\quad\Longleftrightarrow\quad
s\le\frac{\sig^2}{r-1}-\sig(t-2),
\]
and
\[
s+1-\frac{\sig^2}{r+\sig(t-2)-1}\le1
\quad\Longleftrightarrow\quad
s\le\frac{\sig^2}{r+\sig(t-2)-1}.
\]
These are precisely the upper bounds on $s$ in
Definition~\ref{def:admissible-class}\textup{(ii)}--\textup{(iv)},
which proves necessity.

We now prove sufficiency.  Under condition \textup{(i)}, every component
of $H$ belongs to $\mathcal H_+$.  Lemma~\ref{lem:positive-kernel-coronal}
therefore gives $A(H)+\ph I\succeq0$, and
Lemma~\ref{lem:spectral-coronal-new} yields $\lambda_2(G)\le\sig$.

Under one of conditions \textup{(ii)}--\textup{(iv)}, write
$H=H_0\udot R$ as prescribed.  The strict condition on $r$ and
Lemma~\ref{lem:negative-coronal-calculus} show that
$A(H_0)+\ph I$ has exactly one negative eigenvalue and that
$\bbone$ belongs to its range.
Every component of $R$ has a positive semidefinite shifted adjacency
matrix with $\bbone$ in its range by
Lemma~\ref{lem:positive-kernel-coronal}.
Consequently, $n_-(A(H)+\ph I)=1$ and
$\bbone\in\operatorname{range}(A(H)+\ph I)$.
The coronal formulas and the stated bound on $s$ give
$\scoronal(H)=\scoronal(H_0)+s\le1$ by the same equivalences above.
Lemma~\ref{lem:spectral-coronal-new} then yields $\lambda_2(G)\le\sig$.
\end{proof}

\section*{Data availability}
\phantomsection\label{sec:data-availability}
The code and reference output for the finite computational
verifications in this paper are publicly available at
{\def\UrlBreaks{\do\.\do\/\do\-}\url{https://github.com/Chaohui-Chen/golden-section-graphs-verification}}.

\section*{Declaration of generative AI and AI-assisted technologies}
During the preparation of this manuscript, the authors used OpenAI's ChatGPT
for assistance with English-language editing, comparison of manuscript
versions, and drafting supplementary mathematical explanations and calculations. The authors reviewed and revised the
AI-assisted text and take full responsibility for the content of the
manuscript, including all mathematical statements and proofs. No
AI-generated material was treated as a substitute for mathematical proof
or independent verification.

\end{document}